\documentclass[10pt,reqno]{amsart}
\usepackage{mathrsfs}
\usepackage{amssymb,amsmath,amscd,amsthm,amsfonts}
\newcommand{\bea}{\begin{eqnarray}}
\newcommand{\eea}{\end{eqnarray}}
\newcommand{\bna}{\begin{eqnarray*}}
\newcommand{\ena}{\end{eqnarray*}}

\newcommand{\s}{\sigma}

\newtheorem{Theorem}{Theorem}[section]
\newtheorem{Lemma}[Theorem]{Lemma}

\begin{document}

\bigskip
\centerline{\sc \large Solvable Base Change}
\centerline{\sc \large and Rankin-Selberg Convolutions}
\bigskip
\centerline{ Tim Gillespie}

\openup 0.9\jot

\bigskip

\bigskip

\centerline{\scshape Abstract}
\begin{quote}
{\footnotesize In this paper we define a Rankin-Selberg $L$-function
attached to automorphic cuspidal representations of $GL_m({\Bbb
A}_E)\times GL_{m'}({\Bbb A}_{F})$ over solvable algebraic number
fields $E$ and $F$ which are invariant under the Galois action, using a result proved by C.S. Rajan, and prove a prime
number theorem for this $L$-function. 
}
\end{quote}

\bigskip
\section{Introduction}
A prime number theorem for Rankin-Selberg L-functions has already been studied by several authors.  In the classical case of \cite{LiuYe4} the proof requires much of what is known about the classical Rankin-Selberg L-function: meromorphic continuation, functional equation and the zero free region due to Moreno \cite{Mor}.  Unlike the prime number theorem for the zeta function, where non-vanishing on the line $Re(s)=1$ is equivalent to the asymptotic formula 
\begin{equation} \sum_{n\leq x}\Lambda(n)\sim x \nonumber  \end{equation} 
this non-vanishing result does not suffice to give the same estimate for any Rankin-Selberg L-function, and this is the reason for the self-contragredient assumption made in \cite{LiuYe4}.   
 If one uses the automorphic induction functor to define a "Rankin-Selberg product" as in \cite{GillJi} then the main obstacle is to obtain the multiplicity of the poles of such a product.  We now recall the basic notation of \cite{GillJi}, so let $\pi$ be an automorphic cuspidal representation of $GL_n(\mathbb{A}_{E})$ where $E$ is a cyclic algebraic number field of prime degree $\ell$.  Suppose also that $\pi^{\sigma}\cong \pi$ for $\sigma$ a generator of the Galois group $Gal(E/\mathbb{Q})$, then we have the factorization 
 \begin{equation}L(s,\pi)=L(s,\pi_{\mathbb{Q}})L(s,\pi_{\mathbb{Q}}\otimes\eta_{E/\mathbb{Q}})...L(s,\pi_{\mathbb{Q}}\otimes\eta_{E/\mathbb{Q}}^{\ell-1})  \nonumber  \end{equation} 
 where $\{\pi_{\mathbb{Q}}\otimes\eta_{E/\mathbb{Q}}^i\}_{i=0,...,\ell-1}$ are automorphic cuspidal representations on $GL_n(\mathbb{A}_{\mathbb{Q}})$ and $\eta_{E/\mathbb{Q}}$ is an idele class character obtained using the reciprocity isomorphism $$\mathbb{A}_{\mathbb{Q}}^{\times}/\mathbb{Q}^{\times}N_{E/\mathbb{Q}}(\mathbb{A}_E^{\times})\cong Gal(E/\mathbb{Q})$$Similarly, let $\pi'$ be an automorphic cuspidal representation of $GL_m(\mathbb{A}_F)$ where $F$ is a cyclic algebraic number field of prime degree $q$ and $\pi'\cong \pi'^{\tau}$ where $\tau$ is a generator of $Gal(F/\mathbb{Q})$.   Then as before we have 
 \begin{equation} L(s,\pi')=L(s,\pi'_{\mathbb{Q}})L(s,\pi'_{\mathbb{Q}}\otimes\psi_{F/\mathbb{Q}})...L(s,\pi'_{\mathbb{Q}}\otimes\psi_{F/\mathbb{Q}}^{q-1})  \nonumber  \end{equation}
 for an idele class character $\psi_{F/\mathbb{Q}}$.  Then we define the Rankin-Selberg L-function over different fields by 
 \begin{equation} L(s,\pi\times_{BC}\tilde{\pi}')=\prod_{i=0}^{\ell-1}\prod_{j=0}^{q-1}L(s,\pi_{\mathbb{Q}}\otimes\eta_{E/\mathbb{Q}}^{i}\times\widetilde{\pi_{\mathbb{Q}}'\otimes\psi_{F/\mathbb{Q}}^j}) \nonumber \end{equation}
 \begin{equation}=L(s,\boxplus_{i=0}^{\ell-1}(\pi_{\mathbb{Q}}\otimes\eta_{E/\mathbb{Q}}^{i})\times\boxplus_{j=0}^{q-1}(\widetilde{\pi_{\mathbb{Q}}'\otimes\psi_{F/\mathbb{Q}}^{j}}))  \nonumber \end{equation}
\begin{equation}=L(s,AI_{E/\mathbb{Q}}(\pi)\times AI_{F/\mathbb{Q}}(\tilde{\pi}'))     \nonumber \end{equation}  
where $AI_{K/\mathbb{Q}}$ denotes the automorphic induction functor for any number field $K/\mathbb{Q}$.  The Euler product for $L(s,\pi\times_{BC}\pi')$ converges absolutely for $Re(s)>1$, so that we can write 
\begin{equation}  \frac{L'}{L}(s,\pi\times_{BC}\tilde{\pi}')=-\sum_{n\geq 1}\frac{\Lambda(n)a_{\pi\times_{BC}\tilde{\pi}'}(n)}{n^{s}}  \nonumber \end{equation}  
see the next section for a precise definition of $a_{\pi\times_{BC}\pi'}(n)$.  For future reference we let 
\begin{equation} T=\{(\sigma_{\mathbb{Q}},\sigma_{\mathbb{Q}}')|\sigma_{\mathbb{Q}}\in BC_{E/\mathbb{Q}}^{-1}(\pi), \sigma_{\mathbb{Q}}'\in BC_{F/\mathbb{Q}}^{-1}(\pi'),  \sigma_{\mathbb{Q}}\cong \sigma_{\mathbb{Q}}'\otimes|\det|^{i\tau}\exists \tau\in \mathbb{R}\}  \nonumber \end{equation} 
 where $BC_{K/\mathbb{Q}}$ denotes the base change functor for any number field $K$.  Our first result is a prime number theorem for $L(s,\pi\times\tilde{\pi}')$ in the cyclic case when $\ell=q$.  

\begin{Theorem}Let $\pi$ and $\pi'$ be unitary automorphic cuspidal representations of $GL_n(\mathbb{A}_E)$ and $GL_m(\mathbb{A}_F)$, respectively, with $E/\mathbb{Q}$ and $F/\mathbb{Q}$ of prime degree $\ell$ such that $E\neq F$.  Suppose that $\pi$ and $\pi'$ are invariant under the action of $Gal(E/\mathbb{Q})$ and $Gal(F/\mathbb{Q})$, respectively, and with notation as above suppose at least one of $\pi_{\mathbb{Q}}$ or $\pi_{\mathbb{Q}}'$ is self-contragredient.  Then
\begin{equation}  \sum_{n\leq x}\Lambda(n)a_{\pi\times_{BC}\tilde{\pi}'}(n)=    \nonumber \end{equation}
\begin{equation}=\left\{ \begin{array}{l}  \frac{x^{1+i\tau(\pi,\pi')}}{1+i\tau(\pi,\pi')}+O\{x\exp(-c\sqrt{\log x})\} \text{if $T\neq \phi$ and $BC_{EF/E}(\pi)$ is cuspidal}\\
\frac{\ell x^{1+i\tau(\pi,\pi')}}{1+i\tau(\pi,\pi')}+O\{x\exp(-c\sqrt{\log x})\}\text{ if $T\neq \phi$ and $BC_{EF/E}(\pi)$ is not cuspidal} \\
O\{x\exp(-c\sqrt{\log x})\} \text{ if $T=\phi$}  \end{array}  \right.\nonumber \end{equation}     \end{Theorem} 
\medskip 
{\it Remark:}  By Lemma 4.1 of \cite{GillJi} if $T$ is nonempty there is a unique $\tau(\pi,\pi')$ so that $\pi_{\mathbb{Q}}\otimes\eta_{E/\mathbb{Q}}^{i_0}\cong \pi_{\mathbb{Q}}'\otimes\psi_{F/\mathbb{Q}}^{j_0}\otimes|\det|^{i\tau(\pi,\pi')}$ for some $0\leq i_0,j_0\leq \ell-1$.  This follows from multiplicity one for characters and the fact that $\eta_{E/\mathbb{Q}}$ and $\psi_{F/\mathbb{Q}}$ have finite order.   
\medskip 

More generally let $E/\mathbb{Q}$ and $F/\mathbb{Q}$ be any solvable Galois extensions of degrees $\ell$ and $\ell'$, and let $\pi$ and $\pi'$ denote unitary automorphic cuspidal representations of $GL_n(\mathbb{A}_E)$ and $GL_m(\mathbb{A}_F)$, respectively.  Moreover suppose that both $\pi$ and $\pi'$ admit base change from $\mathbb{Q}$; in other words assume that the sets $BC_{E/\mathbb{Q}}^{-1}(\pi)$ and $BC_{F/\mathbb{Q}}^{-1}(\pi')$ are nonempty.  Then by Theorem 2 of \cite{Raj} we can write
 \begin{equation}BC_{E/\mathbb{Q}}^{-1}(\pi)=\{\pi_{\mathbb{Q}}\otimes\chi_i\}_{i\in I}   \nonumber \end{equation}
for some idele class characters $\chi_i$ trivial on $\mathbb{Q}^{\times}N_{E_{ab}/\mathbb{Q}}(\mathbb{A}_{E_{ab}}^{\times})$, where $E_{ab}$ denotes the fixed field of the commutator subgroup $[Gal(E/\mathbb{Q}),Gal(E/\mathbb{Q})]$.  Similarly we can write 
\begin{equation}  BC_{F/\mathbb{Q}}^{-1}(\pi')=\{\pi_{\mathbb{Q}}'\otimes\psi_j\}_{j\in J} \nonumber \end{equation}  
for some idele class characters of $\mathbb{A}_{\mathbb{Q}}^{\times}$ trivial on $\mathbb{Q}^{\times}N_{F_{ab}/\mathbb{Q}}(\mathbb{A}_{F_{ab}}^{\times})$.  Consider the towers of extensions coming from the cylic composition factors of $Gal(E/\mathbb{Q})$ and $Gal(F/\mathbb{Q})$ using the Galois correspondence.
\begin{equation}E=E_0\supset E_1 \supset ...\supset E_{k}\supset E_{k+1}=\mathbb{Q}  
\end{equation}  
\begin{equation} F=F_0\supset F_1 \supset ... \supset F_{r} \supset F_{k+1}=\mathbb{Q}   \end{equation}  
with $[E_i,E_{i+1}]=\ell_{i+1}$ of prime degree for $0\leq i \leq k$ and $[F_j:F_{j+1}]=q_{j+1}$ of prime degree for $0\leq j \leq t$.  We will actually need stronger assumptions on the Galois invariance of the representations in the fibers $BC_{E/E_i}^{-1}(\pi)$ for all $i$.  More precisely suppose that $\pi$ is invariant under the action of $Gal(E/E_1)$ and that the representations in the fiber $BC_{E/E_i}^{-1}(\pi)$ are invariant under the action of $Gal(E_i/E_{i+1})$ for any $2 \leq i \leq k$, then we define the Rankin-Selberg L-function over the fields $E$ and $F$ by 
\begin{eqnarray}  L(s,\pi\times_{BC}\pi')=\prod_{i\in I}\prod_{j\in J}L(s,\pi_{\mathbb{Q}}\otimes\chi_i\times \pi_{\mathbb{Q}}'\otimes \psi_j)   \nonumber \\ =L(s,AI_{E/\mathbb{Q}}(\pi)\times AI_{F/\mathbb{Q}}(\pi)) \nonumber \end{eqnarray}   
 To simplify notation first consider the two step case
  \begin{eqnarray}  E=E_0\supset E_1 \supset E_2=\mathbb{Q}  \end{eqnarray}
Then by assumption the $\ell_1$ distinct representations 
\begin{equation}BC_{E/E_1}^{-1}(\pi)=\{\pi_{E_1}\otimes\eta_{E/E_1}^{i}\}_{i=0}^{\ell_1-1}  \nonumber \end{equation}  
are all invariant under the action of $Gal(E_1/E_2)$ and from the proof of Theorem 2 in \cite{Raj} this forces $\eta_{E/E_1}^{\sigma}\cong \eta_{E/E_1}$ for all $\sigma\in Gal(E_1/E_2)$ so that $\eta_{E/E_1}=\eta_{E/\mathbb{Q}}\circ N_{E_1/\mathbb{Q}}$ for some idele character on $\mathbb{A}_{\mathbb{Q}}^{\times}$ trivial on $\mathbb{Q}^{\times}N_{E_{ab}/\mathbb{Q}}(\mathbb{A}_{E_{ab}}^{\times})$.  Thus we can write 
\begin{equation}BC_{E/\mathbb{Q}}^{-1}=\{\pi_{\mathbb{Q}}\otimes\eta_{E/\mathbb{Q}}^{i}\otimes\eta_{E_1/\mathbb{Q}}^{j}\}_{ {0 \leq j \leq \ell_2-1}\atop{0\leq i \leq \ell_1-1}}  \nonumber \end{equation} 
for some cuspidal automorphic $\pi_{\mathbb{Q}}$ on $GL_n(\mathbb{A}_{\mathbb{Q}})$. Note the above representations are distinct, which can be seen using the fact that $BC_{E/E_1}(\pi_{\mathbb{Q}}\otimes\eta_{E/\mathbb{Q}}^{i})=BC_{E/E_1}(\pi_{\mathbb{Q}})\otimes\eta_{E/E_1}^{i}$. In other words if we have $\pi_{\mathbb{Q}}\otimes\eta_{E/\mathbb{Q}}^{i_1}\otimes\eta_{E_1/\mathbb{Q}}^{j_1}\cong \pi_{\mathbb{Q}}\otimes\eta_{E/\mathbb{Q}}^{i_2}\otimes\eta_{E_1/\mathbb{Q}}^{j_2}$ then the preceding remark implies that $i_1=i_2$ mod($\ell_1$) so that $j_1=j_2$ mod($\ell_2$).  Thus inductively we get that in the general case we have the $\ell$ distinct representations which lift to $\pi$ from $\mathbb{Q}$
\begin{equation} B C_{E/\mathbb{Q}}^{-1}(\pi)= \{\pi_{\mathbb{Q}}\otimes(\otimes_{a=0}^{k}\eta_{E_a/\mathbb{Q}}^{i_a})\}_{i_a=0}^{\ell_{a+1}-1}  \nonumber \end{equation}    
If we make similar Galois invariance assumptions for $\pi'$ then we can write
\begin{equation}BC_{F/\mathbb{Q}}^{-1}(\pi')=\{ \pi_{\mathbb{Q}}'\otimes(\otimes_{b=0}^{t}\psi_{F_b/\mathbb{Q}}^{j_b})\}_{j_b=0}^{q_{b+1}-1} \nonumber \end{equation}
and these are also distinct.  By a similar calculation as in Theorem 1.2 of \cite{GillJi} we put a group structure on the above representations and show that the set of twisted equivalent pairs divides the order of this group.  From this we obtain a prime number theorem for $L(s,\pi\times_{BC}\tilde{\pi}')$ 
\begin{Theorem} Let E and F be solvable Galois extensions of degrees $\ell$ and $\ell'$ with $(\ell, \ell')=1$.  Let $\pi$ and $\pi'$ be unitary automorphic cuspidal representations on $GL_n(\mathbb{A}_{E})$ and $GL_m(\mathbb{A}_F)$ respectively.  Suppose that the representations in the fibers  $BC_{E/E_i}^{-1}(\pi)$, $BC_{F/F_j}^{-1}(\pi')$ are invariant under the action of $Gal(E_i/E_{i+1})$ and $Gal(F_j/F_{j+1})$, respectively for $1\leq i \leq k$, $1\leq j \leq t$. Also suppose that for some $\pi_{\mathbb{Q}}\in BC^{-1}_{E/\mathbb{Q}}(\pi)$ that $\pi_{\mathbb{Q}}$ is self-contragredient, then    
\begin{equation} \sum_{n\leq x}a_{\pi\times_{BC}\tilde{\pi}'}(n)\Lambda(n)=
\end{equation}  
\begin{equation}  \left\{ \begin{array}{l}\frac{x^{1+i\tau(\pi,\pi')}}{1+i\tau(\pi,\pi')}+O\{x\exp(-c\sqrt{\log x})\} \text{  if $T$ is nonempty}    \\  O\{x\exp(-c\sqrt{\log x})\} \text{  if $T$ is empty }   \end{array}\right.  \end{equation}   \end{Theorem}
The method used in proving Theorem 1.1 is to calculate the multiplicity of a pole of $L(s,\pi\times_{BC}\pi')$ and apply the main theorem in \cite{LiuYe4}.  We rely heavily on the description of the fibers of the base change map as proved in \cite{Raj}, and using Theorem 2 from \cite{Raj} combined with Lemma 4.1 from \cite{GillJi} we get that there is at most one distinct pole of $L(s,\pi\times_{BC}\pi')$ with multiplicity equal to one.  

\section{Notation}
We will use the L-functions as in \cite{Jacq}.  For $\pi$ an automorphic cuspidal representation on $GL_n(\mathbb{A}_{E})$ with $E/\mathbb{Q}$ Galois recall that one can define $L(s,\pi)$ as a product of local factors 
\begin{equation}L(s,\pi)=\prod_{p}\prod_{\nu|p}\prod_{i=1}^{n}(1-\alpha_{\pi}(i,\nu)p^{-f_ps})^{-1}  \nonumber  \end{equation} 
where $\{\alpha_\pi(i,\nu)\}_{i=1}^{n}$ is a collection of complex numbers for any place $\nu$ of $E$ and $f_p$ denotes the modular degree of any place $\nu$ lying over $p$.  The above product converges absolutely for $Re(s)>1$ so we can write  
\begin{equation}  \frac{L'}{L}(s,\pi)=-\sum_{n\geq 1}\frac{\Lambda(n)a_{\pi}(n)}{n^s}  \nonumber  \end{equation}
where $\Lambda(n)$ denotes the Von-Mangoldt function, and for $n=p^{f_pk}$ a prime power
\begin{equation}  a_{\pi}(n)=f_p\sum_{\nu|p}\sum_{i=1}^{n}\alpha_{\pi}(i,\nu)^k  \nonumber  \end{equation}
The Galois group $Gal(E/\mathbb{Q})$ acts on $\pi$ by $\pi^{\sigma}(g)=\pi(g^{\sigma^{-1}})$ for $g\in GL_n(\mathbb{A}_{E})$.   In the special case when $E$ is a cyclic algebraic number field of prime degree $\ell$, if we let $<\sigma>=Gal(E/\mathbb{Q})$  by the results in \cite{AC} either $\pi\cong \pi^{\sigma}$ in which case $\pi$ is  the base change lift of exactly $\ell$ non-equivalent cuspidal representations $\{\pi_{\mathbb{Q}}\otimes\eta_{E/\mathbb{Q}}^j\}_{j=0}^{\ell-1}$ where $\pi_{\mathbb{Q}}\otimes\eta_{E/\mathbb{Q}}^j$ is on $GL_n(\mathbb{A}_{\mathbb{Q}})$,  or we have $\pi\ncong\pi^{\sigma}$ so that
\begin{equation}L(s,\pi)=L(s,\pi_{\mathbb{Q}})  \nonumber  \end{equation}
with $\pi_{\mathbb{Q}}$ an automorphic cuspidal representation of $GL_{n\ell}(\mathbb{A}_{\mathbb{Q}})$.     
We will also need some results regarding Rankin-Selberg L-functions.  We will use the Rankin-Selberg $L$-functions $L(s, \pi \times
\widetilde\pi')$ as developed by Jacquet, Piatetski-Shapiro, and
Shalika \cite{JacPiaShal}, Shahidi \cite{Sha1}, and Moeglin and
Waldspurger \cite{MoeWal}, where $\pi$ and $\pi '$ are unitary automorphic cuspidal representations of $GL_m(\mathbb{A}_E)$ and
$GL_{m'}(\mathbb{A}_E)$, respectively.  Recall $L(s,\pi\times\tilde{\pi}')$ is defined as the product of local factors
\begin{equation}\prod_{p}L_p(s,\pi\times\tilde{\pi}')=\prod_p\prod_{\nu|p}\prod_{i=1}^{m}\prod_{j=1}^{m'}\Big(1-\alpha_{\pi}(i,\nu)\overline{\alpha_{\pi'}(j,\nu)}p^{-f_ps}\Big)^{-1}  \nonumber \end{equation} 
where $f_p$ denotes the modular degree of any place $\nu|p$.  We will need the following properties of the $L$-functions
$L(s,\pi\times\widetilde{\pi}')$.

\medskip

{\bf RS1}. The Euler product defining $L(s,\pi\times\tilde{\pi}')$ converges absolutely for $\s>1$ (Jacquet and Shalika
\cite{JacSha1}).

\medskip

{\bf RS2}. $L(s,\pi\times\tilde{\pi}')$ admits meromorphic continuation to all of $\mathbb{C}$, and if we denote $\alpha(g)=|\det(g)|$ and  $\pi'\not\cong
\pi\otimes\alpha^{it}$ for any $t\in{\Bbb R}$, then
$L(s,\pi\times\widetilde{\pi}')$ is holomorphic. When $m=m'$ and
$\pi' \cong \pi\otimes \alpha^{i\tau_0} $ for some $\tau_0\in\Bbb
R$, the only poles of $L(s, \pi \times \widetilde\pi ')$ are simple
poles at $s=i\tau_0$ and $1+i\tau_0$ (Jacquet and Shalika \cite{JacSha1},
Moeglin and Waldspurger \cite{MoeWal}).

\medskip

{\bf RS3}.
$L(s,\pi\times\widetilde{\pi}')$ is non-zero in $\s\ge 1$ (Shahidi
\cite{Sha1}).  Furthermore, if at least one of $\pi$ or $\pi'$ is
self-contragredient, it is zero-free in the region
\begin{equation}
\sigma > 1-\frac{c}{\log(Q_{\pi}Q_{\pi'}(|t|+2))}, \quad |t|\geq 1
\end{equation}
where $c$ is an explicit constant depending only on $m$ and $n$ (see
Sarnak \cite{Sa}, Moreno \cite{Mor} or Gelbert, Lapid and Sarnak
\cite{GLS}).
  
\section{Proof of Theorem 1.1}  
{\it Proof. } First note that since $E$ and $F$ are of prime degree, if $E\neq F$ then we have an isomorphism of Galois groups $Gal(EF/\mathbb{Q})\cong Gal(E/\mathbb{Q})\times Gal(F/\mathbb{Q})$ given by restriction.  We also have that $EF/\mathbb{Q}$ is a solvable extension of degree $\ell^2$, so that the base change map $BC_{EF/\mathbb{Q}}$
exists, and we denote $\pi_{EF}=BC_{EF/E}(\pi)$.  Suppose first that $\pi_{EF}$ is cuspidal, then from \cite{Raj} we get $BC_{EF/\mathbb{Q}}^{-1}(\pi_{EF})=\{\pi_{\mathbb{Q}}\otimes\chi_{i}\}_{i\in I}$ for some idele class characters of the quotient  
\begin{equation}  \chi_{i}:\mathbb{A}_{\mathbb{Q}}^{\times}/(\mathbb{Q}^{\times}N_{EF/\mathbb{Q}}(\mathbb{A}_{EF}^{\times}))\longrightarrow \mathbb{C}^{\times}  \nonumber \end{equation}   

Note that there are $\ell^2$ distinct representations which lift to $\pi_{EF}$ if $\pi_{EF}$ is cuspidal.  To see this take 
\begin{equation}  \nonumber  BC_{EF/E}^{-1}(\pi_{EF})=\{\pi\otimes\eta_{EF/E}^j\}_{j=0}^{\ell}
\end{equation}  
and these representations are distinct.  Now let $C_{K}=\mathbb{A}_{K}^{\times}/K^{\times}$ for any number field $K$, and consider the character $\psi_{F/\mathbb{Q}}\circ N_{E/\mathbb{Q}}$ which is a character on $C_{E}$ trivial on $N_{EF/E}(C_{EF})$.  Thus we have that for some $0\leq i\leq \ell-1$, $\eta_{EF/E}^{i}=\psi_{F/\mathbb{Q}}\circ N_{E/\mathbb{Q}}$ and if the character on the right hand side is trivial we get that $N_{E/\mathbb{Q}}(C_{E})\subseteq ker(\psi_{F/\mathbb{Q}})=N_{F/\mathbb{Q}}(C_F)$ so by class field theory $F\subseteq E$ which is a contradiction, so that $\psi_{F/\mathbb{Q}}\circ N_{E/\mathbb{Q}}$ is nontrivial and so has order $\ell$.   Hence $\eta_{EF/E}=\eta_{EF/\mathbb{Q}}\circ N_{E/\mathbb{Q}}$ for some class field theory character $\eta_{EF/\mathbb{Q}}$, in other words $\eta_{EF/E}$ lies in the image of the base change map so that we can take 
\begin{equation} BC_{E/\mathbb{Q}}^{-1}(\pi\otimes\eta_{EF/E}^j)=\{\pi^{j}\otimes\eta_{E/\mathbb{Q}}^{i}\}_{i=0}^{\ell-1}   \nonumber  \end{equation}   
and again these are distinct for each fixed $j$.  Consider the collection $\{\pi^{j}\otimes\eta_{E/\mathbb{Q}}^{i} \}_{0\leq i,j\leq \ell-1}$, which all lift to $\pi_{EF}$ by the transitivity of the base change map; moreover they are distinct since given $\pi^{j_1}\otimes\eta_{E/\mathbb{Q}}^{i_1}\cong \pi^{j_2}\otimes\eta_{E/\mathbb{Q}}^{i_2}$ then applying $BC_{E/\mathbb{Q}}$ gives 
\begin{equation} \pi\otimes\eta_{EF/E}^{j_1}\cong\pi\otimes\eta_{EF/E}^{j_2}       \nonumber  \end{equation}      
which implies $j_1=j_2$ mod($\ell$) so that $i_1=i_2$  mod($\ell$) and  this proves the claim.      
Finally using the isomorphism of Galois groups we get that any character on $\mathbb{A}_{\mathbb{Q}}^{\times}/(\mathbb{Q}^{\times}N_{EF/\mathbb{Q}}(\mathbb{A}_{EF}^{\times}))$ may be written as $\eta_{E/\mathbb{Q}}^{i}\otimes\psi_{F/\mathbb{Q}}^{j}$ for some $0\leq i,j\leq \ell-1$. From this and the preceding remarks it follows that the representations 
\begin{equation}\{ \pi_{\mathbb{Q}}\otimes\eta_{E/\mathbb{Q}}^i\otimes\psi_{F/\mathbb{Q}}^j\}_{0\leq i,j\leq \ell-1 }  \nonumber  \end{equation}  
are distinct.  Now suppose the set $T$ of twisted equivalent pairs is nonempty, so that for some $0\leq i_0,j_0\leq \ell-1$ and $\tau_0\in\mathbb{R}$ we have
\begin{equation} \pi_{\mathbb{Q}}\otimes\eta_{E/\mathbb{Q}}^{i_0}\cong\pi'_{\mathbb{Q}}\otimes\psi_{F/\mathbb{Q}}^{j_0}\otimes|\det|^{i\tau_0}    \nonumber   \end{equation}   
If we have another twisted equivalent pair, say 
\begin{equation} \pi_{\mathbb{Q}}\otimes\eta_{E/\mathbb{Q}}^{i_1}\cong \pi_{\mathbb{Q}}'\otimes\psi_{F/\mathbb{Q}}^{j_1}\otimes|\det|^{i\tau_1}   \nonumber \end{equation}   
Then by Lemma 1 of \cite{GillJi} we may suppose that $\tau_0=\tau_1$, thus we get 
\begin{eqnarray} \pi_{\mathbb{Q}}\otimes\eta_{E/\mathbb{Q}}^{i_1}\otimes\psi_{F/\mathbb{Q}}^{j_0} \cong\pi_{\mathbb{Q}}'\otimes\psi_{F/\mathbb{Q}}^{j_1}\otimes\psi_{F/\mathbb{Q}}^{j_0}\otimes|\det|^{i\tau_0}   \nonumber   \\ 
\cong \pi_{\mathbb{Q}}\otimes\eta_{E/\mathbb{Q}}^{i_0}\otimes\psi_{F/\mathbb{Q}}^{j_1}  \nonumber  \end{eqnarray}
so that $i_1=i_0$ mod($\ell$) and $j_1=j_0$ mod($\ell$) as desired.  Finally suppose that $\pi_{EF}$ is not cuspidal, then from \cite{AC} we get that $\ell|n$ and 
\begin{equation} \pi\otimes\eta_{EF/E}\cong \pi  \nonumber    \end{equation} 
for some nontrivial character of $\mathbb{A}_{E}^{\times}/(E^{\times}N_{EF/E}(\mathbb{A}_{EF}^{\times}))$.  As before we get that $\eta_{EF/E}=\psi_{F/\mathbb{Q}}^{k}\circ N_{E/\mathbb{Q}}$ for some $1 \leq r \leq \ell-1$, hence $\pi_{\mathbb{Q}}\cong \pi_{\mathbb{Q}}\otimes\eta_{E/\mathbb{Q}}^{s}\otimes\psi_{F/\mathbb{Q}}^{r}$
for some $0\leq s \leq \ell-1$ again from \cite{AC}.  As before suppose that  
\begin{equation}\pi_{\mathbb{Q}}\otimes\eta_{E/\mathbb{Q}}^{i_0}\cong \pi_{\mathbb{Q}}'\otimes\psi_{F/\mathbb{Q}}^{j_0}\otimes|\det|^{i\tau_0}  \end{equation}   
Then by a simple calculation we get another twisted equivalent pair 
\begin{equation} \pi_{\mathbb{Q}}\otimes\eta_{E/\mathbb{Q}}^{s+i_0}\cong\pi_{\mathbb{Q}}'\otimes\psi_{F/\mathbb{Q}}^{j_0-r}\otimes|\det|^{i\tau_0}   \nonumber \end{equation}  
and this pair is distinct from (3.7), so by Lemma 4.2 from \cite{GillJi} we get $|T|=\ell$.  The rest of the proof follows as in the proof of Theorem 1.2 of \cite{GillJi}.  $\square$  
\smallskip
\section{Proof of Theorem 1.2}  
For completeness we first state a Lemma which is the analogue of Lemma 4.1 of \cite{GillJi}, we omit the proof as it is almost identical to that given in \cite{GillJi}  
\begin{Lemma}Suppose that $\pi_{\mathbb{Q}}\otimes(\otimes_{a=0}^{k}\eta_{E_a/\mathbb{Q}}^{i_{a,0}})\cong\pi_{\mathbb{Q}}'\otimes(\otimes_{b=0}^{t}\psi_{F_b/\mathbb{Q}}^{j_{b,0}})\otimes|\det|^{i\tau_0}$ for some $0 \leq i_{a,0} \leq \ell_{a+1}$ and $0 \leq j_{b,0}\leq q_{b+1}$ then 
\begin{equation} \pi_{\mathbb{Q}}\otimes(\otimes_{a=0}^{k}\eta_{E_a/\mathbb{Q}}^{i_{a,0}})\cong \pi_{\mathbb{Q}}'\otimes(\otimes_{b=0}^{t}\psi_{F_b/\mathbb{Q}}^{j_b})|\det|^{i\tau}   \nonumber \end{equation} 
implies that $j_b=j_{b,0}$ for $b=0,...,t$ and $\tau=\tau_0$.  Moreover, if for some $i_a$ and $j_b$ with $a=0,...,k$, $b=0,...,t$
\begin{equation}  
\pi_{\mathbb{Q}}\otimes(\otimes_{a=0}^{k}\eta_{E_a/\mathbb{Q}}^{i_a})\cong\pi_{\mathbb{Q}}'\otimes(\otimes_{b=0}^{t}\psi_{F_b/\mathbb{Q}}^{j_b})\otimes|\det|^{i\tau}   \nonumber \end{equation}  
then $\tau=\tau_0$
\nonumber \end{Lemma} 
By Lemma 4.1 if  the set $T$ is nonempty the exponent $\tau$ is uniquely determined and we denote this by $\tau_{(\pi,\pi')}$.  We will use the notation as in the introduction.  Since the representations 
\begin{equation}BC_{E/\mathbb{Q}}^{-1}(\pi)=\{\pi_{\mathbb{Q}}\otimes(\otimes_{a=0}^{k}\eta_{E_a/\mathbb{Q}}^{i_a})\}_{i_a=0}^{\ell_{a+1}-1}  \nonumber \end{equation}  
are distinct, we get a well-defined group operation by setting 
\begin{equation} (\pi_{\mathbb{Q}}\otimes(\otimes_{a=0}^{k}\eta_{E_a/\mathbb{Q}}^{i_a}))*(\pi_{\mathbb{Q}}\otimes(\otimes_{a=0}^{k}\eta_{E_a/\mathbb{Q}}^{i_a'}))   \nonumber \end{equation}  
 \begin{equation}=\pi_{\mathbb{Q}}\otimes(\otimes_{a=0}^{k}\eta_{E_a/\mathbb{Q}}^{i_a+i_a'})  \nonumber \end{equation}  
and we denote this group of order $\ell$ by $(\mathcal{G},*)$. Now suppose the set $T$ is nonempty, then  
\begin{equation}\sigma_{\mathbb{Q}}\cong \sigma_{\mathbb{Q}}'\otimes |\det|^{i\tau_{(\pi,\pi')}}  \nonumber \end{equation}
 for any $(\sigma_{\mathbb{Q}}, \sigma_{\mathbb{Q}}')\in T$, and by relabeling if necessary we may assume that $(\pi_{\mathbb{Q}},\pi_{\mathbb{Q}}'\otimes(\otimes_{b=0}^{t}\psi_{F_b/\mathbb{Q}}^{j_{b,0}}))\in T$ for some $\pi_{\mathbb{Q}}'\otimes(\otimes_{b=0}^{t}\psi_{F_b/\mathbb{Q}}^{j_{b,0}}))\in BC_{F/\mathbb{Q}}^{-1}(\pi')$.  Moreover given two twisted equivalent pairs 
 \begin{eqnarray}\pi_{\mathbb{Q}}\otimes(\otimes_{a=0}^{k}\eta_{E_a/\mathbb{Q}}^{i_a})\cong \pi_{\mathbb{Q}}'\otimes(\otimes_{b=0}^{t}\psi_{F_b/\mathbb{Q}}^{j_b})\otimes|\det|^{i\tau_{(\pi,\pi')}} \nonumber \\
 \pi_{\mathbb{Q}}\otimes(\otimes_{a=0}^{k}\eta_{E_a/\mathbb{Q}}^{i_a'})\cong \pi_{\mathbb{Q}}'\otimes(\otimes_{b=0}^{t}\psi_{F_b/\mathbb{Q}}^{j_b'})\otimes |\det|^{i\tau_{(\pi,\pi')}}\end{eqnarray}   
we obtain 
\begin{equation}\pi_{\mathbb{Q}}\otimes(\otimes_{a=0}^{k}\eta_{E_a/\mathbb{Q}}^{i_a-i_a'})\cong \pi_{\mathbb{Q}}'\otimes(\otimes_{b=0}^{t}\psi_{F_b/\mathbb{Q}}^{j_b})\otimes(\otimes_{a=0}^{k}\eta_{E_a/\mathbb{Q}}^{-i_a'})\otimes|\det|^{i\tau_{(\pi,\pi')}}\nonumber \end{equation}
\begin{equation}\cong \pi_{\mathbb{Q}}\otimes(\otimes_{b=0}^{t}\psi_{F_b/\mathbb{Q}}^{j_b-j_b'})\cong\pi_{\mathbb{Q}}'\otimes(\otimes_{b=0}^{t}\psi_{F_b/\mathbb{Q}}^{j_{b,0}+j_b-j_b'})\otimes|\det|^{i\tau_{(\pi,\pi')}}  \nonumber \end{equation}
so it follows that the subset $\mathcal{H}\subset \mathcal{G}$ defined by 
\begin{equation}  \mathcal{H}=\{\sigma_{\mathbb{Q}}\in \mathcal{G}| (\sigma_{\mathbb{Q}},\sigma_{\mathbb{Q}}')\in T, \exists \sigma_{\mathbb{Q}}'\in BC_{F/\mathbb{Q}}^{-1}(\pi')\}  \nonumber \end{equation}      
 forms a subgroup of $\mathcal{G}$ so by LaGrange's theorem we get $|\mathcal{H}|$ divides $\ell$.  Moreover Lemma 4.1 gives that for $(\sigma_{\mathbb{Q}},\sigma_{\mathbb{Q}}')\in T$ then $\sigma_{\mathbb{Q}}$ is twisted equivalent to at most one $\sigma_{\mathbb{Q}}'\in BC_{F/\mathbb{Q}}^{-1}(\pi')$, hence $|T|=|\mathcal{H}|$.  Thus applying the same argument above to the collection 
 \begin{equation} \{\pi_{\mathbb{Q}}'\otimes(\otimes_{b=0}^{t}\psi_{F_b/\mathbb{Q}}^{j_b})\}_{j_b=0}^{q_{b+1}-1}    \nonumber\end{equation}    
we get that the cardinality of $T$ also divides $\ell'$, so that $|T|=1$ since $(\ell,\ell')=1$.  
Finally, assuming $\pi_{\mathbb{Q}}$ to be self-contragredient and using the equality 
 \begin{equation}L(s,\pi_{\mathbb{Q}}\times\chi\times\widetilde{\pi_{\mathbb{Q}}'\otimes\xi})=L(s,\pi_{\mathbb{Q}}\times\widetilde{\pi_{\mathbb{Q}}'\otimes\chi^{-1}\xi}) \nonumber \end{equation} 
valid for any finite order idele class characters we can apply the zero-free region to obtain the desired error term, and the theorem follows.  $\square$

\section{  Acknowledgements }  The author would like to thank Professors Yang Bo Ye and Muthu Krishnamurthy for their helpful suggestions and encouragement.  The author would also like to thank Professor Wang Song at the Beijing Istitute for his helpful advice, and the University of Iowa for it support.   
{100}
\end{document}